\documentclass[numreferences]{kluwer}
\usepackage{amsthm, amssymb, amsfonts}
\begin{document}
\begin{article}
\newtheorem{prop}{Proposition}
\begin{opening}
\title{Semi-classical twists for $\mathfrak{sl}_{3}$ and 
$\mathfrak{sl}_{4}$ boundary $r-$matrices of Cremmer-Gervais type}
\author{M. \surname{Samsonov }\email{samsonov@pink.phys.spbu.ru}}
\institute{Theoretical Department\\
Institute of Physics\\
St. Petersburg State University\\
198904, St. Petersburg\\
Russia}
\runningtitle{Semi-classical twists for $\mathfrak{sl}_{3}$ and 
$\mathfrak{sl}_{4}$ boundary $r-$matrices}
\runningauthor{M. Samsonov}
\begin{ao}
Theoretical Department\\
Institute of Physics\\
198904, St. Petersburg\\
Russia
\end{ao}
\begin{abstract}
We obtain explicit formulas for the semi-classical twists deforming
the coalgebraic structure of $U(\mathfrak{sl}_{3})$ and $U(\mathfrak{sl}_{4})$. 
In rank $2$ and $3$ the corresponding universal $R-$matrices 
quantize the boundary $r-$matrices of Cremmer-Gervais type 
defining Lie Frobenius structures on the maximal parabolic subalgebras in 
$\mathfrak{sl}_{n}$.  
\end{abstract}
\classification{Mathematics Subject Classifications(2005)}{16W30, 17B37, 81R50}
\end{opening}
\section{Introduction}
An interesting subclass of the boundary $r-$matrices is given by the generalized 
Jordanian $r-$matrices of Cremmer-Gervais type which are the boundary points 
of $SL_{n}$ adjoint action orbits containing skew-symmetric Cremmer-Gervais $r-$matrices 
\cite{EH,GG}. Explicitly, these boundary 
points in the $\mathfrak{sl}_{n}$ Cartan-Weyl basis are the following   
\begin{equation}
r_{\mathfrak{p}}=
\sum_{p=1}^{n-1}D_{p}\wedge E_{p,p+1}+\sum_{i<j}\sum_{m=1}^{j-i-1}
E_{i,j-m+1}\wedge E_{j,i+m},
\label{boundary}
\end{equation}
where 
$$
D_{p}=
\displaystyle\frac{n-p}{n}(E_{11}+E_{22}+\cdots+E_{pp})-\displaystyle\frac pn (E_{p+1,p+1}+E_{p+2,p+2}+\cdots+E_{nn}).\\[2ex]
$$
Each $r_{\mathfrak{p}}$ defines the structure of Lie Frobenius algebra 
(a Lie algebra with a nondegenerate $2-$coboundary \cite{GG}) on the 
maximal parabolic subalgebra $\mathfrak{p}\subset\mathfrak{sl}_{n}$ 
generated by the Cartan subalgebra and all the simple root generators 
excluding $E_{n,n-1}$. In this letter we construct the twists 
and the universal $R_{\mathfrak{p}}-$matrices if $n=3,4$.
The concrete $R-$matrices $R^{V}$ arise when one restricts  
$R_{\mathfrak{p}}={F_{\mathfrak{p}}}_{21}F_{\mathfrak{p}}^{-1}$ 
to a particular $\mathfrak{sl}_{n}$ representation $V$. It was observed 
\cite{EH} that if $V_{n}$ is the linear space of polynomials 
of the degree $\le n$ then $R^{V_{n}}$ are related to rational 
degeneration of the Ueno-Shibukawa operators and one obtains 
$R^{V_{n}}$ explicitly. Thus this letter gives an answer to the next 
question what the universal $R-$matrices are in two first cases of 
physical interest, the genuine Cremmer-Gervais case $n=3$ and quantization of the 
complexified conformal algebra $o(4,2)_{\mathbb{C}}\approx \mathfrak{sl}_{4}$. 
The paths of constructing the parabolic twists, as we name 
$F_{\mathfrak{p}}$ following \cite{Sam}, are almost parallel in both 
cases and below we describe them uniformly. We consider 
$U^{\cal K}_{q}(\hat{\mathfrak{sl}}_{n-1})$, the 
Drinfeld-Jimbo quantization $U_{q}(\hat{\mathfrak{sl}}_{n-1})$ deformed 
by the abelian twist ${\cal K}$, and propose a method of constructing 
the affine twists allowing nontrivial specialization in the 
limit $q\rightarrow 1$. The construction we follow is based on 
factorization of a chosen singular trivial twist 
$F_{n-1}:=(W_{n-1}\otimes W_{n-1})
\Delta_{\cal K}(W^{-1}_{n-1})$,
with 
$W_{n-1}\in U_{q}^{\cal K}(\hat{\mathfrak{sl}}_{n-1})$, into singular 
$\Phi_{W_{n-1}}^{-1}$ and nonsingular $F_{n-1}^{\rm aff}$ 
parts such that $F_{n-1}=\Phi_{W_{n-1}}^{-1}\cdot F_{n-1}^{\rm aff}$. 
Among all the factorizations we find the one such that 
$F_{n-1}^{\rm aff}$ is a twist equivalent to an affine version of 
the Cremmer-Gervais twist or its $\mathfrak{sl}_{4}$ analog \cite{IO,KST,KM}. 
The final step in our approach is rational degeneration of 
$F^{\rm aff}_{n-1}$ which we denote $\overline{F^{\rm aff}_{n-1}}$ and 
construction  of a homomorpism $\iota_{n-1}$ such that 
$(\iota_{n-1}\otimes\iota_{n-1})(\overline{F^{\rm aff}_{n-1}})$ is 
a twist on $U^{\Psi_{n}}(\mathfrak{sl}_{n})$, where $\Psi_{n}$ 
turns out to be the semi-classical twist found in \cite{KLM}. 
The final twist on $U(\mathfrak{sl}_{n})$ is obtained as the composition 
$(\iota_{n-1}\otimes\iota_{n-1})(\overline{F^{\rm aff}_{n-1}})\cdot\Psi_{n}$.  
\section{Quantum affine twists for $U_{q}(\hat{\mathfrak{sl}_{2}})$ and 
$U_{q}(\hat{\mathfrak{sl}_{3}})$} Fix the central charge $c=1$ in the 
defining relations of $U_{q}(\hat{\mathfrak{sl}}_{n})$. Let 
$A=\left(\frac{2(\alpha_{i}|\alpha_{j})}{(\alpha_{i}|\alpha_{i})}\right)_{i,j=1}^{n}$ be the 
Cartan matrix of $\mathfrak{sl}_{n}$ and $\{\alpha_{i}\}_{i=0,\ldots, n-1}$ 
be the set of all simple roots of $\hat{\mathfrak{sl}}_{n}$  with the symmetric 
scalar product $(\cdot|\cdot)$ on it. With these assumptions the relations 
defining $U_{q}(\hat{\mathfrak{sl}}_{n})$ over the field of rational functions 
$\mathbb{Q}(q)$ are the following 
$$
\begin{array}{lcr}
q^{h_{\alpha_{0}}}q^{h_{\alpha_{1}}}\ldots q^{h_{\alpha_{n-1}}}=1,
&& q^{h_{\alpha_{i}}}q^{h_{\alpha_{j}}}=q^{h_{\alpha_{j}}}q^{h_{\alpha_{i}}}\\[2ex]
q^{h_{\alpha_{i}}}e_{\pm\alpha_{j}}q^{-h_{\alpha_{i}}}=
q^{\pm(\alpha_{i}|\alpha_{j})}e_{\pm\alpha_{j}}, && 
[e_{\alpha_{i}},e_{-\alpha_{j}}]=\delta_{ij}
\displaystyle\frac{q^{h_{{\alpha}_{i}}}-q^{-h_{{\alpha}_{i}}}}{q-q^{-1}}\\[2ex]
e_{\pm\alpha_{i}}e_{\pm\alpha_{j}}
=e_{\pm\alpha_{j}}e_{\pm\alpha_{i}}&\mbox{if}&|i-j|\nequiv 1\pmod{n}\\[2ex]
[e_{\pm\alpha_{i}},[e_{\pm\alpha_{i}},e_{\pm\alpha_{j}}]_{q}]_{q^{-1}}=0
&\mbox{if}& |i-j|\equiv 1\pmod{n}, n>2\\[2ex]
[e_{\pm\alpha_{i}},[e_{\pm\alpha_{i}},[e_{\pm\alpha_{i}},e_{\pm\alpha_{j}}
]_{q^{2}}]]_{q^{-2}}=0 &\mbox{if}& |i-j|\equiv 1\pmod{n}, n=2
\end{array}
$$
and the $q-$commutator is defined as usual 
$$
\begin{array}{lcl}
[C,D]_{q}=CD-qDC,&& [C,D]\equiv[C,D]_{1}.
\end{array}
$$
$U_{q}(\hat{\mathfrak{sl}}_{n})$ is a Hopf algebra. The comultiplication
is uniquely defined by fixing its value 
on $q-$Chevalley generators
$$
\begin{array}{lclclcl}
\Delta(e_{{\alpha}_{i}})&=&q^{-h_{\alpha_{i}}}\otimes e_{{\alpha}_{i}}
+e_{{\alpha}_{i}}\otimes 1, &&
\Delta(e_{{-\alpha}_{i}})&=&e_{{-\alpha}_{i}}\otimes q^{h_{\alpha_{i}}}+1\otimes e_{{-\alpha}_{i}}\\[2ex]
\end{array}
$$
$$
\Delta(q^{h_{{\alpha}_{i}}})=q^{h_{{\alpha}_{i}}}\otimes q^{h_{{\alpha}_{i}}}.
$$  
It is convenient to use the following notation
$$
(x)_{(1)}:=x\otimes 1,\hphantom{aaa} (x)_{(2)}:=~1\otimes x
$$
and its modified version 
$$
\begin{array}{c}
(x)_{<1>}:=({\rm id}\otimes{\rm pr}_{K_{n-1}})\Delta(x),\hphantom{aaa} 
(x)_{<2>}:=({\rm pr}_{K_{n-1}}\otimes{\rm id})\Delta(x)\\[2ex] 
(x)_{<3>}:=(({\rm id}-{\rm pr}_{K_{n-1}})\otimes ({\rm id}-{\rm pr}_{K_{n-1}}))
\Delta(x)
\end{array}
$$
where ${\rm pr}_{K_{n-1}}$ the 
$U_{q}(\hat{\mathfrak{sl}}_{n})$ projector to 
$K_{n-1}:=\mathbb{Q}(q)[q^{\pm h_{\alpha_{i}}}]_{i=1,\ldots,n-1}$.

\subsection{An affine twist for $U_{q}(\hat\mathfrak{sl}_{2})$}
Introduce a topological Hopf algebra ${\cal D}^{(2)}[[\zeta]]$ 
as a completion in the formal series topology of the following 
subalgebra   
$$
{\cal D}^{(2)}:=\{\sum_{l_{1},l_{2},l_{3}\ge 0}c_{l_{1},l_{2},l_{3}}
(e_{-\alpha})^{l_{1}}q^{\pm l_{2} h_{\alpha}}(e_{\delta-\alpha})^{l_{3}}|
c_{l_{1},l_{2},l_{3}}\in\mathbb{Q}(q)\}
\subset U_{q}(\hat{\mathfrak{sl}}_{2})
$$
and consider a trivial twist $F_{2}$:
$$
\begin{array}{c} 
F_{2}=(W_{2}\otimes W_{2})\Delta(W^{-1}_{2}),\\[2ex] 
W_{2}=\exp_{q^{2}}(\displaystyle\frac{\zeta}{1-q^{2}}\mathop{}e_{\delta-\alpha})\exp_{q^{-2}}
(-\displaystyle\frac{\zeta^{2}}{1-q^{2}}\mathop{}q^{-h_{\alpha}}e_{-\alpha})\in{\cal D}^{(2)}[[\zeta]];\\[2ex]
\exp_{q^{2}}(z):=\displaystyle\sum_{k\ge 0}\frac{z^{k}}{(k)_{q^{2}}!},\hphantom{a}(k)_{q^{2}}!=\frac{1-q^{2}}{1-q^{2}}\cdots\frac{1-q^{2n}}{1-q^{2}}
\end{array}
$$
where $\delta=\alpha_{0}+\alpha$. 
Along with the nonsingular $q-$exponent we make use of its singular version    
$
{\rm e}_{q^{2}}(z):=~\exp_{q^{2}}(\frac{1}{1-q^{2}}\mathop{}z)
$ 
and of several results of $q-$calculus \cite{FK, Kac} concerning its 
properties:
\begin{itemize} 
\item If $[x,y]_{q^{2}}=0$ then 
\begin{equation}
\begin{array}{lcl}
{\rm e}_{q^{2}}(x+y)={\rm e}_{q^{2}}(y)\cdot{\rm e}_{q^{2}}(x),&&
({\rm e}_{q^{2}}(x))^{-1}={\rm e}_{q^{-2}}(q^{-2}\mathop{}x).
\label{mult} 
\end{array}
\end{equation}
\item If $[u,[u,v]]_{q^{2}}=[v,[u,v]]_{q^{-2}}=0$ then
\begin{equation}
{\rm e}_{q^{2}}(u)\cdot{\rm e}_{q^{2}}(v)={\rm e}_{q^{2}}(v)\cdot
{\rm e}_{q^{2}}(\displaystyle\frac{1}{1-q^{2}}\mathop{}[u,v])\cdot{\rm e}_{q^{2}}(u).
\label{5terms}
\end{equation}
\item The Heine's formula holds
\begin{equation}
\displaystyle(1-u)_{q^{2}}^{(v)}=
{\rm e}_{q^{2}}(u)({\rm e}_{q^{2}}(uq^{-2v}))^{-1}
\label{Heine}
\end{equation}
$$
\displaystyle(1-u)_{q^{2}}^{(v)}
:=1+\sum_{k>0}\frac{(-v)_{q^{2}}(-v+1)_{q^{2}}\cdots (-v+k-1)_{q^{2}}}{(k)_{q^{2}}!}u^{k}
$$ 
where $(-v+l)_{q^{2}}=(q^{-2v+2l}-1)/(q^{2}-1)$.
\end{itemize}
\begin{prop}
\label{prop1}
$$
F_{2}=\Phi_{W_{2}}^{-1}\cdot F_{2}^{\rm aff}
$$
where $F_{2}^{\rm aff}\in (K_{1}\otimes{\cal D}^{(2)})[[\zeta]]$
is a twist and 
$$
\Phi_{W_{2}}=
{\rm Ad} (W_{2}\otimes W_{2})({\rm e}_{q^{2}}(-\zeta^{3}
\mathop{}e_{\delta-\alpha}q^{-h_{\alpha}}\otimes q^{-h_{\alpha}}e_{-\alpha}))=
{\rm Ad} (W_{2}\otimes W_{2})(\Phi_{2}).
$$
\end{prop}
\begin{pf}
Using (\ref{mult}) and the coproducts
$$
\begin{array}{lcl}
\Delta(e_{\delta-\alpha})=q^{h_{\alpha}}\otimes e_{\delta-\alpha}+e_{\delta-\alpha}\otimes 1,
&&
\Delta(e_{-\alpha})=e_{-\alpha}\otimes q^{h_{\alpha}}+1\otimes e_{-\alpha}.
\end{array}
$$
$F_{2}$ factorizes into the product of $q-$exponents    
\begin{equation}
\begin{array}{r}
(W_{2}\otimes W_{2})\cdot
{\rm e}_{q^{2}}(\zeta^{2}\mathop{}(q^{-h_{\alpha}}e_{-\alpha})_{<1>})\cdot
{\rm e}_{q^{2}}(\zeta^{2}\mathop{}(q^{-h_{\alpha}}e_{-\alpha})_{<2>})\cdot\\[2ex]
{\rm e}_{q^{-2}}(q^{-2}\zeta\mathop{}(e_{\delta-\alpha})_{<1>})\cdot
{\rm e}_{q^{-2}}(q^{-2}\zeta\mathop{}(e_{\delta-\alpha})_{<2>}).
\end{array}
\label{triv}
\end{equation}
Let us focus on the following product  
\begin{equation}
{\rm e}_{q^{2}}(\zeta^{2}\mathop{}(q^{-h_{\alpha}}e_{-\alpha})_{<2>})\cdot
{\rm e}_{q^{-2}}(q^{-2}\zeta\mathop{}(e_{\delta-\alpha})_{<1>}).
\label{switch}
\end{equation}
Applying (\ref{5terms}) to 
$$
\begin{array}{l}
{\rm e}_{q^{-2}}(q^{-2}\zeta\mathop{}(e_{\delta-\alpha})_{<1>})\cdot
{\rm e}_{q^{-2}}(q^{-2}\zeta^{2}\mathop{}(q^{-h_{\alpha}}e_{-\alpha})_{<2>})\\[2ex]
={\rm e}_{q^{-2}}(q^{-2}\zeta^{2}\mathop{}(q^{-h_{\alpha}}e_{-\alpha})_{<2>})\cdot
\Phi^{-1}_{2}\cdot{\rm e}_{q^{-2}}
(q^{-2}\zeta\mathop{}(e_{\delta-\alpha})_{<1>})
\end{array}
$$
we can flip the $q-$exponents in (\ref{triv}) and after collecting them 
using (\ref{mult}) we come to  
\begin{equation}
F_{2}=\Phi^{-1}_{W_{2}}\cdot{\rm e}_{q^{2}}(U)\cdot
{\rm e}_{q^{-2}}(q^{-2}\mathop{}q^{(h_{\alpha})_{(1)}}U)
\label{U}
\end{equation}
where
$
U=\zeta\mathop{}(e_{\delta-\alpha})_{(2)}+
\zeta^{2}\mathop{}(q^{-h_{\alpha}}e_{-\alpha})_{<2>}.
$\\
The final step in the proof is application of the Heine's formula
 that leads to an explicit form of $F_{2}^{\rm aff}$
$$
F_{2}^{\rm aff}=(1-\zeta\mathop{}(e_{\delta-\alpha})_{(2)}-
\zeta^{2}\mathop{}(q^{-h_{\alpha}}e_{-\alpha})_{<2>})_{q^{2}}
^{(-\frac 12 h_{\alpha})_{(1)}}
$$
that satisfies the 
Drinfeld equation according to \cite{KST}.
\end{pf} 
\newproof{rem}{Remark}
\begin{rem}
By Proposition \ref{prop1} $F_{\mathfrak{p}}^{\rm aff}$ is equivalent 
to $\Phi_{2}$ and the latter is the Cremmer-Gervais twist 
as it is seen from \cite{KM}. In Proposition \ref{prop2} we 
give a proof for the case $U_{q}(\hat{\mathfrak{sl}_{3}})$ and 
indicate what are the simplifications one needs to take 
into consideration to adopt it for $U_{q}(\hat{\mathfrak{sl}_{2}})$.
\end{rem}
\subsection{An affine twist for $U_{q}(\hat{\mathfrak{sl}}_{3})$}
Let us deform the coalgebraic structure of  
$U_{q}(\hat{\mathfrak{sl}}_{3})$ by the following 
abelian twist 
$$
{\cal K}=q^{\frac 49 h_{\alpha}\otimes h_{\alpha}+
\frac 29 h_{\alpha}\otimes h_{\beta}+
\frac 59 h_{\beta}\otimes h_{\alpha}+\frac 79 h_{\beta}\otimes h_{\beta}}
$$
the convenience of this choice will be justified by Proposition \ref{prop2} 
(see also \cite{IO})
$$
\Delta_{\cal K}(\hat{e}_{-\beta})=q^{h_{\alpha+\beta}^{\perp}}
\otimes\hat{e}_{-\beta}+\hat{e}_{-\beta}\otimes 1,\hphantom{a}
\Delta_{\cal K}(\hat{e}_{-\alpha})=
q^{-h_{\beta}^{\perp}}\otimes\hat{e}_{-\alpha}+\hat{e}_{-\alpha}\otimes 1
$$
$$
\begin{array}{c}
\Delta_{\cal K}(\hat{e}_{\alpha})=\hat{e}_{\alpha}
\otimes q^{h_{\alpha+\beta}^{\perp}}+1\otimes\hat{e}_{\alpha}\\[2ex]
\Delta_{\cal K}(\hat{e}_{\delta-\alpha-\beta})
=\hat{e}_{\delta-\alpha-\beta}\otimes q^{-h_{\beta}^{\perp}}
+1\otimes\hat{e}_{\delta-\alpha-\beta}
\end{array}
$$
$$
\begin{array}{c}
\Delta_{\cal K}(\hat{e}_{\delta-\beta})
=\hat{e}_{\delta-\beta}\otimes q^{-h_{\alpha}^{\perp}}
+1\otimes\hat{e}_{\delta-\beta}+(1-q^{2})\mathop{}
\hat{e}_{\alpha}\otimes\hat{e}_{\delta-\alpha-\beta}
q^{h_{\alpha+\beta}^{\perp}}\\[2ex]
\Delta_{\cal K}(\hat{e}_{-\alpha-\beta})
=q^{-h_{\alpha}^{\perp}}\otimes\hat{e}_{-\alpha-\beta}
+\hat{e}_{-\alpha-\beta}\otimes 1+
(1-q^{-2})\mathop{}\hat{e}_{-\beta}q^{-h_{\beta}^{\perp}}\otimes\hat{e}_{-\alpha}
\end{array}
$$
where $\delta=\alpha_{0}+\alpha+\beta$
$$
\begin{array}{llll}
h_{\alpha}^{\perp}=\frac 23 h_{\alpha}+\frac 43 h_{\beta},&
h_{\beta}^{\perp}=\frac 43 h_{\alpha}+\frac 23 h_{\beta},&
h_{\alpha+\beta}^{\perp}=h_{\beta}^{\perp}-h_{\alpha}^{\perp},&
\hat{e}_{-\beta}=q^{-\frac 12\mathop{}h_{\beta}^{\perp}}e_{-\beta}
\end{array}
$$
$$
\begin{array}{lll}
\hat{e}_{-\alpha}=q^{\frac 12\mathop{}h_{\beta}^{\perp}}e_{-\alpha},&
\hat{e}_{\alpha}=e_{\alpha}q^{\frac 12\mathop{}h_{\alpha+\beta}^{\perp}},
&
\hat{e}_{\delta-\alpha-\beta}=q^{-\frac 12\mathop{}h_{\alpha+\beta}^{\perp}}
e_{\delta-\alpha-\beta}
\end{array}
$$
$$
\begin{array}{ll}
\hat{e}_{\delta-\beta}=\hat{e}_{\alpha}\hat{e}_{\delta-\alpha-\beta}-
q^{2}\hat{e}_{\delta-\alpha-\beta}\hat{e}_{\alpha},&
\hat{e}_{-\alpha-\beta}=\hat{e}_{-\beta}\hat{e}_{-\alpha}-q^{-2}
\hat{e}_{-\alpha}\hat{e}_{-\beta}.
\end{array}
$$
\begin{rem}
${\cal K}$ preserves the composite root generators 
in the sense that  
$$
\begin{array}{lcl}
\hat{e}_{\delta-\alpha}=e_{\alpha}e_{\delta-\alpha-\beta}-qe_{\delta-\alpha-\beta}
e_{\alpha},&&
\hat{e}_{-\alpha-\beta}=e_{-\beta}e_{-\alpha}-q^{-1}e_{-\alpha}e_{-\beta}.
\end{array}
$$ 
\end{rem}
Denote by $K_{2}^{\prime}$ the minimal algebra containing 
$K_{2}\cup\{\hat{e}_{\alpha}\}$ over $\mathbb{Q}(q)$.  
Introduce a subalgebra ${\cal D}^{(3)}\subset U_{q}(\hat{\mathfrak{sl}}_{3})$ 
generated by all finite linear combinations with coefficients in $\mathbb{Q}(q)$:
$$
\sum_{l_{1},\ldots,l_{5}\ge 0} c_{l_{1},\ldots, l_{5}} (\hat{e}_{-\beta})^{l_{1}}
(\hat{e}_{-\alpha-\beta})^{l_{2}}(\hat{e}_{-\alpha})^{l_{3}}K_{2}^{\prime}
(\hat{e}_{\delta-\beta})^{l_{4}}(\hat{e}_{\delta-\alpha-\beta})^{l_{5}}.
$$
If we complete ${\cal D}^{(3)}$ up to 
${\cal D}^{(3)}[[\zeta]]$ and choose the following element   
$$
\begin{array}{l}
W_{3}=\\
{\rm e}_{q^{2}}(\zeta\mathop{}q^{h_{\alpha}^{\perp}}
\hat{e}_{\delta-\beta})
{\rm e}_{q^{2}}(-\displaystyle\frac{q\zeta^{2}}{1-q^{2}}\mathop{}q^{h_{\beta}^{\perp}}\hat{e}_{\delta-\alpha-\beta})
({\rm e}_{q^{2}}(-\zeta^{2}\mathop{}\hat{e}_{-\beta})
{\rm e}_{q^{2}}(-\zeta\mathop{}\hat{e}_{-\alpha})
{\rm e}_{q^{2}}(-\zeta^{2}\mathop{}\hat{e}_{-\beta}))^{-1}
\end{array}
$$
defining the trivial twist 
$$
F_{3}=(W_{3}\otimes W_{3})\Delta_{\cal K}(W^{-1}_{3}) 
$$
then similarly to Proposition \ref{prop1} we can formulate 
\begin{prop}
\label{prop2}
\begin{equation}
F_{3}=\Phi_{W_{3}}^{-1}\cdot F_{3}^{\rm aff}
\label{claim}
\end{equation}
where 
$F_{3}^{\rm aff}\in (K_{2}^{\prime}\otimes{\cal D}^{(3)})[[\zeta]]$ is a twist 
and $\Phi_{W_{3}}={\rm Ad}(W_{3}\otimes W_{3})(\Lambda_{1}\Lambda_{2}\Lambda_{3})$
$$
\begin{array}{l}
\Lambda_{1}=\exp_{q^{2}}(-\displaystyle\frac{q^{-1}\zeta^{4}}{(1-q^{2})^{2}}
\mathop{}\hat{e}_{\delta-\alpha-\beta}\otimes\hat{e}_{-\beta})\cdot
\exp_{q^{2}}(-\displaystyle\frac{q^{-1}\zeta^{3}}{(1-q^{2})^{2}}\mathop{}\hat{e}_{\delta-\alpha-\beta}\otimes \hat{e}_{-\alpha})\\[2ex]
\Lambda_{2}=
\exp_{q^{2}}(\displaystyle\frac{\zeta^{4}}{(1-q^{2})^{2}}\mathop{}\hat{e}_{\delta-\beta}^{\prime}\otimes\hat{e}_{-\alpha-\beta}),\hphantom{a}\hat{e}_{\delta-\beta}^{\prime}=
[\hat{e}_{\alpha},\hat{e}_{\delta-\alpha-\beta}]\\[2ex]
\Lambda_{3}=
\exp_{q^{2}}(-\zeta(q-q^{-1})\mathop{}\hat{e}_{\alpha}\otimes \hat{e}_{-\beta}).
\end{array}
$$ 
\end{prop}
\begin{pf}
A general idea of factorization of $F_{3}$ is to move factors from 
$\Delta_{\cal K}(W^{-1}_{3})$ containing identity in the second tensor 
factor to the left in order to form $\Phi_{W_{3}}^{-1}$.
Using explicit form of the coproducts and the commutation relations
$[\hat{e}_{\alpha},\hat{e}_{\delta-\beta}]=[\hat{e}_{\delta-\alpha-\beta},\hat{e}_{\delta-\beta}]=0$,
we expand $\Delta_{\cal K}(W^{-1}_{3})$ into the product of 
$q-$exponents using (\ref{mult}). According to the strategy of 
factorization we follow, first we flip the following two $q-$exponents 
from $\Delta_{\cal K}(W^{-1}_{3})$ by introducing the factor 
$\Lambda_{3}$:
$$
\begin{array}{l}
{\rm e}_{q^{2}}(-\zeta^{2}\mathop{}(\hat{e}_{-\beta})_{<2>})\cdot
{\rm e}_{q^{2}}(-\zeta\mathop{}(\hat{e}_{-\alpha})_{<1>})=\\[2ex]
{\rm Ad}(\Lambda_{3})^{-1}
\left({\rm e}_{q^{2}}(-\zeta\mathop{}(\hat{e}_{-\alpha})_{<1>})\cdot
{\rm e}_{q^{2}}(-\zeta^{2}\mathop{}q^{-h_{\beta}^{\perp}}\otimes
\hat{e}_{-\beta})\right).
\end{array}
$$
The latter is seen from the relation
$$ 
\zeta^{2}\mathop{}q^{h_{\alpha+\beta}^{\perp}}\otimes \hat{e}_{-\beta}+
\zeta\mathop{}\hat{e}_{-\alpha}\otimes 1=
{\rm Ad}(\Lambda_{3})^{-1}
\left(\zeta\mathop{}\hat{e}_{-\alpha}\otimes 1+\zeta^{2}\mathop{}q^{-h_{\beta}^{\perp}}
\otimes\hat{e}_{-\beta}\right)
$$
if one notices
$$
[\hat{e}_{\alpha},\hat{e}_{-\alpha}]=\displaystyle
\frac{q^{h_{\alpha+\beta}^{\perp}}-
q^{-h_{\beta}^{\perp}}}{q-q^{-1}}
$$
and applies (\ref{mult}). Then we move $\Lambda_{3}$ to the 
right of $\Delta_{\cal K}(W^{-1}_{3})$ by applying (\ref{5terms}) and 
the commutation relations 
$$
\begin{array}{rcrcr}
[\hat{e}_{\delta-\beta},\hat{e}_{\alpha}]=0, &&
[\hat{e}_{\delta-\beta},\hat{e}_{-\beta}]=0,&& [\hat{e}_{\delta-\beta},\hat{e}_{\delta-\alpha-\beta}]=0\\[2ex]
\end{array}
$$
until we arrive at the following form 
of $F_{3}$  
$$
\begin{array}{r}
(W_{3}\otimes W_{3})\Lambda_{3}^{-1}
({\rm e}_{q^{2}}(-\zeta^{2}\mathop{}\hat{e}_{-\beta})\cdot
{\rm e}_{q^{2}}(-\zeta\mathop{}\hat{e}_{-\alpha})\cdot
{\rm e}_{q^{2}}(-\zeta^{2}\mathop{}\hat{e}_{-\beta}))_{<1>}\cdot\hphantom{aaaaaaaaaaaaaaa}\\[2ex]
{\rm e}_{q^{2}}(-\zeta^{2}\mathop{}q^{-h_{\beta}^{\perp}}\otimes
\hat{e}_{-\beta})\cdot
{\rm e}_{q^{2}}(-\zeta\mathop{}(\hat{e}_{-\alpha})_{<2>})\cdot
\exp_{q^{2}}(-q^{-1}\zeta^{2}\mathop{}\hat{e}_{\alpha}q^{-h_{\beta}^{\perp}}\otimes 
\hat{e}_{-\alpha-\beta})\cdot\hphantom{aaaaa}\\[2ex]
{\rm e}_{q^{-2}}(-\displaystyle\frac{q^{-1}\zeta^{2}}{1-q^{2}}\mathop{}
(q^{h_{\beta}^{\perp}}
\hat{e}_{\delta-\alpha-\beta})_{<1>})\cdot\Lambda_{3}\cdot
{\rm e}_{q^{-2}}(q^{-2}\zeta\mathop{}(q^{h_{\alpha}^{\perp}}\hat{e}_{\delta-\beta})_{<1>})\cdot\hphantom{aaaaa}\\[2ex]
{\rm e}_{q^{2}}(-\zeta^{2}\mathop{}(\hat{e}_{-\beta})_{<2>})\cdot
{\rm e}_{q^{-2}}(-\displaystyle\frac{q^{-1}\zeta^{2}}{1-q^{2}}\mathop{}(q^{h_{\beta}^{\perp}}
\hat{e}_{\delta-\alpha-\beta})_{<2>})\cdot\hphantom{aaaaa}\\[2ex]
{\rm e}_{q^{-2}}(q^{-2}\zeta\mathop{}(q^{h_{\alpha}^{\perp}}\hat{e}_{\delta-\beta})_{<3>})\cdot
{\rm e}_{q^{-2}}(q^{-2}\zeta\mathop{}(q^{h_{\alpha}^{\perp}}\hat{e}_{\delta-\beta})_{<2>}).\hphantom{aaaaaa}
\end{array}
$$
Next step is the appearance of $\Lambda_{2}^{-1}$ from the relation
$$
\begin{array}{l}
\exp_{q^{2}}(-q^{-1}\zeta^{2}\mathop{}\hat{e}_{\alpha}q^{-h_{\beta}^{\perp}}
\otimes\hat{e}_{-\alpha-\beta})\cdot
{\rm e}_{q^{-2}}(-\displaystyle\frac{q^{-1}\zeta^{2}}{1-q^{2}}\mathop{}(q^{h_{\beta}^{\perp}}
\hat{e}_{\delta-\alpha-\beta})_{<1>})=\\[2ex]
\Lambda_{2}^{-1}\cdot{\rm e}_{q^{-2}}(-\displaystyle\frac{q^{-1}\zeta^{2}}{1-q^{2}}\mathop{}(q^{h_{\beta}^{\perp}}
\hat{e}_{\delta-\alpha-\beta})_{<1>})\cdot
\exp_{q^{2}}(-q^{-1}\zeta^{2}\mathop{}\hat{e}_{\alpha}q^{-h_{\beta}^{\perp}}
\otimes\hat{e}_{-\alpha-\beta})
\end{array}
$$
holding by $[\hat{e}_{\delta-\alpha-\beta},\hat{e}_{\delta-\beta}^{\prime}]_{q^{-2}}=
[\hat{e}_{\alpha},\hat{e}_{\delta-\beta}^{\prime}]_{q^{2}}=0$ and by (\ref{5terms}).
We move $\Lambda_{2}^{-1}$ to the left of $F_{3}$ using 
(\ref{5terms}) and the necessary for it relations  
$$
\begin{array}{lcl}
[\hat{e}_{-\alpha},\hat{e}_{-\alpha-\beta}]=[\hat{e}_{-\beta},\hat{e}_{-\alpha-\beta}]=0,
&&[\hat{e}_{\delta-\beta}^{\prime},\hat{e}_{-\alpha}]=-q\mathop{}
\hat{e}_{\delta-\alpha-\beta}q^{-h_{\beta}^{\perp}}
\end{array}
$$
and form $\Phi_{W_{3}}^{-1}$ which leads to appearance of 
$\Lambda_{1}^{-1}$ and results in the following 
form of 
\begin{equation}
F_{3}=\label{lines}
\end{equation} 
$$
\begin{array}{l}
\Phi^{-1}_{W_{3}} (1\otimes W_{3})\cdot
{\rm e}_{q^{2}}(-\zeta^{2}\mathop{}q^{-h_{\beta}^{\perp}}\otimes\hat{e}_{-\beta})
\cdot{\rm e}_{q^{2}}(-\zeta\mathop{}(\hat{e}_{-\alpha})_{<2>})\cdot\\[2ex]
\exp_{q^{2}}(-q^{-1}\zeta^{2}\mathop{}
\hat{e}_{\alpha}q^{-h_{\beta}^{\perp}}\otimes\hat{e}_{-\alpha-\beta})\cdot\Lambda_{3}\cdot
{\rm e}_{q^{2}}(-\zeta^{2}\mathop{}(\hat{e}_{-\beta})_{<2>})\cdot\\[2ex]
\hphantom{aaa}{\rm e}_{q^{-2}}(-
\displaystyle\frac{q^{-1}\zeta^{2}}{1-q^{2}}\mathop{}
(q^{h_{\beta}^{\perp}}\hat{e}_{\delta-\alpha-\beta})_{<2>})\cdot
{\rm e}_{q^{-2}}(q^{-2}\zeta\mathop{}
(q^{h_{\alpha}^{\perp}}\hat{e}_{\delta-\beta})_{<3>})\cdot\\[2ex]
\hphantom{aaaaaaaaaaaaaaaaaaaaaaaaaaaaaaa}{\rm e}_{q^{-2}}(q^{-2}\zeta\mathop{}
(q^{h_{\alpha}^{\perp}}\hat{e}_{\delta-\beta})_{<2>}).
\end{array}
$$
The last steps will be to note that 
$$
\begin{array}{l}
[\left({\rm e}_{q^{-2}}(-q^{-2}\zeta\mathop{}\hat{e}_{-\alpha})\cdot
{\rm e}_{q^{-2}}(-q^{-2}\zeta^{2}\mathop{}\hat{e}_{-\beta})\right)_{(2)},\\[2ex]
\hphantom{aaaaaaaaaaaaaaa}{\rm e}_{q^{2}}
(-\zeta^{2}\mathop{}q^{-h_{\beta}^{\perp}}\otimes \hat{e}_{-\beta})\cdot
{\rm e}_{q^{2}}(-\zeta\mathop{}(\hat{e}_{-\alpha})_{<2>})]=0
\end{array}
$$
and rewrite $F_{3}$ in the following form with the help of (\ref{5terms}):
$$
\begin{array}{l}
\Phi_{W_{3}}^{-1}\cdot\left(
{\rm e}_{q^{2}}(\zeta\mathop{} q^{h_{\alpha}^{\perp}}\hat{e}_{\delta-\beta})\cdot
{\rm e}_{q^{-2}}(-q^{-2}\zeta^{2}\mathop{} \hat{e}_{-\beta})\right)_{(2)}\cdot\\[2ex]
\hphantom{aa}{\rm e}_{q^{2}}(-\zeta^{2}\mathop{}
q^{-h_{\beta}^{\perp}}\otimes \hat{e}_{-\beta})\cdot
{\rm e}_{q^{-2}}(q^{-2}\zeta\mathop{}q^{h_{\beta}^{\perp}}\otimes
q^{h_{\alpha}^{\perp}}\hat{e}_{\delta-\beta})\cdot\\[2ex]
\hphantom{aaaaa}{\rm e}_{q^{2}}(\zeta\mathop{}q^{h_{\beta}^{\perp}}\otimes
q^{h_{\alpha}^{\perp}}\hat{e}_{\delta-\beta})\cdot
{\rm e}_{q^{2}}
(-\displaystyle\frac{q\zeta^{2}}{1-q^{2}}\mathop{}(q^{h_{\beta}^{\perp}}\hat{e}_{\delta-\alpha-\beta})_{(2)})\cdot
{\rm e}_{q^{2}}(-\zeta\mathop{}
(\hat{e}_{-\alpha})_{<2>})\cdot\\[2ex]
\hphantom{aaaaaaaaaaaaaaaaaaaaaaaa}{\rm e}_{q^{-2}}(-q^{-2}\zeta\mathop{}(\hat{e}_{-\alpha})_{(2)})\cdot
{\rm e}_{q^{-2}}(-q^{-2}\zeta\mathop{}(\hat{e}_{-\beta})_{(2)})\cdots
\end{array}
$$
where we have inserted 
\begin{equation}
{\rm e}_{q^{-2}}(q^{-2}\zeta^{2}\mathop{}q^{h_{\beta}^{\perp}}\otimes
q^{h_{\alpha}^{\perp}}\hat{e}_{\delta-\beta})\cdot
{\rm e}_{q^{2}}(\zeta^{2}\mathop{}q^{h_{\beta}^{\perp}}\otimes
q^{h_{\alpha}^{\perp}}\hat{e}_{\delta-\beta})=1
\label{identity}
\end{equation}
and $\cdots$ means third, forth and fifth lines in (\ref{lines}).
The final transformation is to move the second $q-$exponent 
in (\ref{identity}) 
to the right of $F_{3}$ so that to apply the Heine's formula by the 
same trick as we did for $F_{2}$. $F_{3}=$
$$
\begin{array}{r}
\Phi_{W_{3}}^{-1}\cdot(1-\zeta\mathop{}(q^{h_{\alpha}^{\perp}}\hat{e}_{\delta-\beta})_{(2)}
+\zeta^{2}\mathop{}q^{-(h_{\beta}^{\perp})_{(1)}}(\hat{e}_{-\beta})_{(2)}
)_{q^{2}}^{(-\frac 12\mathop{}h_{\beta}^{\perp})_{(1)}}\cdot
(1+q^{-2}\zeta\mathop{}(\hat{e}_{\alpha})_{(2)})_{q^{-2}}
^{(-\frac 12 h_{\beta}^{\perp})_{(1)}}\cdot\\[2ex]
\exp_{q^{2}}(-q^{-1}\zeta^{2}\mathop{}
\hat{e}_{\alpha}q^{-h_{\beta}^{\perp}}\otimes\hat{e}_{-\alpha-\beta})\cdot
\Lambda_{3}\cdot
{\rm e}_{q^{-2}}(q^{-2}\zeta\mathop{} (q^{h_{\alpha}^{\perp}}\hat{e}_{\delta-\beta})_{<3>})\cdot\\[2ex]
(1-q^{-2}\zeta\mathop{}
(q^{h_{\alpha}^{\perp}}\hat{e}_{\delta-\beta})_{<2>}
+q^{-2}\zeta^{2}\mathop{}(\hat{e}_{\beta})_{(2)})_{q^{-2}}
^{(\frac 12 h_{\alpha+\beta}^{\perp})_{(1)}}.
\end{array}
$$ 
Once we have factored $F_{3}$ we can give a simple proof of the 
Drinfeld equation 
$$
(F_{3}^{\rm aff})_{12}(\Delta\otimes{\rm id})(F_{3}^{\rm aff})=
(F_{3}^{\rm aff})_{23}({\rm id}\otimes\Delta)(F_{3}^{\rm aff}).
$$
From the factorization we know that 
\begin{equation}
(F_{3}^{\rm aff})=(W_{3}\otimes W_{3})\Phi_{3}\Delta_{\cal K}(W^{-1}_{3})
\label{eq}
\end{equation}
where $\Phi_{3}:=\Lambda_{1}\Lambda_{2}\Lambda_{3}$. 
Let us consider the Drinfeld associator
\begin{equation}
{\rm Assoc}(F_{3}^{\rm aff})\equiv(F_{3}^{\rm aff})_{23}({\rm id}\otimes\Delta)(F_{3}^{\rm aff})
((F_{3}^{\rm aff})_{12}(\Delta\otimes{\rm id})(F_{3}^{\rm aff}))^{-1}
\label{assoc2}
\end{equation}
then
\begin{equation}
{\rm Assoc}(F_{3}^{\rm aff})\in (K_{2}^{\prime}\otimes {\cal D}^{(3)}
\otimes{\cal D}^{(3)})[[\zeta]].
\label{assoc3}
\end{equation}
On the other hand by (\ref{eq}) we have 
\begin{equation}
{\rm Assoc}(F_{3}^{\rm aff})=(W_{3}\otimes W_{3}\otimes W_{3}){\rm Assoc}(\Phi_{3})(W^{-1}_{3}
\otimes W^{-1}_{3}\otimes W^{-1}_{3})
\label{assoc4}
\end{equation}
Denote by ${\rm pr}_{K_{2}^{\prime}}$ the projection of 
${\cal D}^{(3)}[[\zeta]]$ to $K_{2}^{\prime}[[\zeta]]$.
Thus by (\ref{assoc3}) 
$${\rm Assoc}(F_{3}^{\rm aff})=
({\rm pr}_{K_{2}^{\prime}}\otimes{\rm id}\otimes {\rm id})
({\rm Assoc}(F_{3}^{\rm aff}))
$$
and  using explicit form of $\Phi_{3}$ 
and $W_{3}$ we deduce that ${\rm Assoc}(F_{3}^{\rm aff})$ is equal to
$$
({\rm pr}_{K_{2}^{\prime}}\otimes{\rm id}\otimes {\rm id})
\left(W_{3}^{\otimes 3}(\Phi_{3})_{23}({\rm id}\otimes \Delta)(\Lambda_{3})
((\Lambda_{3})_{12}(\Lambda_{1})_{23}(\Lambda_{2})_{23}(\Delta\otimes{\rm id})(\Lambda_{3}))^{-1}W_{3}^{-\otimes 3}\right)	  
$$
and the latter is $1\otimes 1\otimes 1$.
\end{pf}
\begin{rem}
In the case of $F_{2}^{\rm aff}$ the proof of the Drinfeld equation 
is similar and one needs to use projection to $K_{1}$ and follow the same lines 
as we did for $F_{3}^{\rm aff}$
\end{rem}
\section{Rational degeneration of $F_{n-1}^{\rm aff},$ $n=3,4$}
\subsection{Rational degeneration of $F_{2}^{\rm aff}$} 
\label{aff22}
Introduce ${\cal D}^{(2)}_{\cal A}[[\zeta]]$ as a completion 
in the formal series topology of the following Hopf subalgebra in 
$U_{q}(\hat{\mathfrak{sl}}_{2})$ 
$$
{\cal D}^{(2)}_{\cal A}:=\{\sum_{l_{1},l_{2},l_{3}\ge 0}c_{l_{1},l_{2},l_{3}}
(e_{-\alpha})^{l_{1}}\left(\frac{q^{\pm h_{\alpha}}-1}{q-1}
\right)^{l_{2}}(e_{\delta-\alpha})^{l_{3}}|
c_{l_{1},l_{2},l_{3}}\in{\cal A}\}
$$
where ${\cal A}:=\mathbb{Q}[q,q^{-1}]_{(q-1)}$ be the ring 
$\mathbb{Q}[q,q^{-1}]$ localized at $(q-1)$ (the subring of rational functions 
nonsingular at $q=1$). Then consider the subalgebra 
${\cal F}^{(2)}_{\cal A}[\zeta]$ generated over ${\cal A}[\zeta]$ by 
the following elements  
$$
\begin{array}{lclcl}
H_{\pm\alpha}=\displaystyle\frac{q^{\pm h_{\alpha}}-1}{q-1},&& 
f_{0}=(q^{-2}-1)\mathop{}q^{-h_{\alpha}}e_{-\alpha},&&
f_{1}=e_{\delta-\alpha}+\zeta\mathop{}q^{-h_{\alpha}}e_{-\alpha}
\end{array}
$$ 
where $q^{-h_{\alpha}}=
1+(q-1)\left(\displaystyle\frac{q^{-h_{\alpha}}-1}{q-1}\right)$.\\ 
The specialization 
${\cal F}^{(2)}_{{\cal A},\mathop{}q=1}[\zeta]=
{\cal F}^{(2)}_{\cal A}[\zeta]/(q-1){\cal F}^{(2)}_{\cal A}[\zeta]$ 
is a Hopf algebra with the following structure  
\begin{equation}
\begin{array}{ccccc}
[\overline{H^{\vphantom{f}}}_{\alpha},
\overline{f^{\vphantom{f}}}_{1}]=-2\overline{f^{\vphantom{f}}}_{1},&& 
[\overline{H^{\vphantom{f}}}_{\alpha},\overline{f^{\vphantom{f}}}_{0}]=
-2\overline{f^{\vphantom{f}}}_{0},&& 
\overline{f^{\vphantom{f}}}_{1}\overline{f^{\vphantom{f}}}_{0}
-\overline{f^{\vphantom{f}}}_{0}\overline{f^{\vphantom{f}}}_{1}=-\zeta
\mathop{}\overline{f^{\vphantom{f}}}_{0}^{2}
\label{relations1}
\end{array}
\end{equation}
\begin{equation}
\Delta(\overline{f^{\vphantom{f}}}_{1})=\overline{f^{\vphantom{f}}}_{1}\otimes 1+1
\otimes\overline{f^{\vphantom{f}}}_{1}
+\overline{H^{\vphantom{f}}}_{\alpha}\otimes\overline{f^{\vphantom{f}}}_{0},
\label{corelations}
\end{equation}
and $\overline{H^{\vphantom{f}}}_{-\alpha}=-\overline{H^{\vphantom{f}}}_{-\alpha}$, 
$\overline{f^{\vphantom{f}}}_{0}$ are primitive. The relations (\ref{relations1}) define
${\cal F}^{(2)}_{{\cal A},\mathop{}q=1}[\zeta]$ as an associative algebra 
with the basis $\{\overline{f^{\vphantom{f}}}_{0}^{k}
\overline{H^{\vphantom{f}}}_{\alpha}^{l}\overline{f^{\vphantom{f}}}_{1}^{m}\}_{k,l,m\ge 0}$ 
over $\mathbb{Q}[\zeta]$ obtained by specialization $q~=~1$ from 
its quantum version. Comparing the semi-classical basis with its quantum 
analog  $\{(f_{0})^{k}
(H_{\alpha})^{l}f_{1}^{m}\}_{k,l,m\ge 0}$ we deduce that $(\ref{relations1})$ are the 
only relations defining  ${\cal F}^{(2)}_{{\cal A},\mathop{}q=1}[\zeta]$.
Completing  ${\cal F}^{(2)}_{{\cal A},\mathop{}q=1}[\zeta]$ up to a topological 
Hopf algebra 
${\cal F}^{(2)}_{{\cal A},\mathop{}q=1}[[\zeta]]:=({\cal F}^{(2)}_{{\cal A},\mathop{}q=1}[\zeta])[[\zeta]]$
we see that $F_{2}^{\rm aff}$ given by
$$
F_{2}^{\rm aff}=(1\otimes 1- \zeta\mathop{}1\otimes f_{1}-\zeta^{2}\mathop{}(h_{\alpha}/2)_{q^{-2}}\otimes 
f_{0})_{q^{2}}
^{\left (-\frac 12 h_{\alpha}\otimes 1\right)}
$$
can be specialized to a twist of ${\cal F}^{(2)}_{{\cal A},\mathop{}q=1}[[\zeta]]$
where it is given by the formula of \cite{KST} 
\begin{equation}
\overline{F_{2}^{\rm aff}}=(1\otimes 1-\zeta\mathop{}
1\otimes \overline{f^{\vphantom{f}}}_{1}-\zeta^{2}\mathop{}
(\overline{H^{\vphantom{f}}}_{\alpha}/2)\otimes 
\overline{f^{\vphantom{f}}}_{0})^{-\frac 12 \overline{H^{\vphantom{f}}}_{\alpha}\otimes 1} .
\label{contracted}
\end{equation}
\subsection{Rational degeneration of $F_{3}^{\rm aff}$}
Let
$$
{\cal D}^{(3)}_{\cal A}:=
\sum_{l_{1},\ldots,l_{5}\ge 0} c_{l_{1},\ldots, l_{5}} (\hat{e}_{-\beta})^{l_{1}}
(\hat{e}_{-\alpha-\beta})^{l_{2}}(\hat{e}_{-\alpha})^{l_{3}}K_{2}^{\prime\prime}
(\hat{e}_{\delta-\beta})^{l_{4}}(\hat{e}_{\delta-\alpha-\beta})^{l_{5}}
$$
where $c_{l_{1},\ldots, l_{5}}\in{\cal A}$. 
$K_{2}^{\prime\prime}$ is an algebra generated by 
$H_{\pm\alpha}^{\perp}=
\displaystyle\frac{q^{\pm h_{\alpha}^{\perp}}-1}{q-1}$, $H_{\pm\beta}^{\perp}=
\displaystyle\frac{q^{\pm h_{\beta}^{\perp}}-1}{q-1}$ and $\hat{e}_{\alpha}$ over
${\cal A}$. Complete ${\cal D}^{(3)}_{\cal A}$ up to 
${\cal D}^{(3)}_{\cal A}[[\zeta]]$ and consider the Hopf subalgebra 
\begin{equation}
{\cal F}^{(3)}_{\cal A}[\zeta]:= 
\{\sum_{l_{1},\ldots,l_{5}\ge 0} c_{l_{1},\ldots, l_{5}} (f_{2})^{l_{1}}
(f_{0})^{l_{2}}(\hat{e}_{-\alpha})^{l_{3}}K_{2}^{\prime\prime}
(f_{3})^{l_{4}}(f_{1})^{l_{5}}|c_{l_{1},\ldots,l_{5}}\in{\cal A}[\zeta]\}
\label{basis321}
\end{equation}
where
$$
\begin{array}{rcl}
f_{0}=(q-q^{-1})\mathop{}\hat{e}_{-\alpha-\beta},&& 
f_{1}=q^{h_{\beta}^{\perp}}\hat{e}_{\delta-\alpha-\beta}
+q^{-1}\zeta\mathop{}\hat{e}_{-\alpha-\beta}
\\[2ex]
f_{2}=(1-q^{-2})\mathop{}\hat{e}_{-\beta},&& 
f_{3}=q^{h_{\alpha}^{\perp}}\hat{e}_{\delta-\beta}-
\zeta\mathop{}\hat{e}_{-\beta}.\\[2ex]
\end{array}
$$
\begin{prop}
$F_{2}^{\rm aff}$ restricts to a twist of  
${\cal F}^{(3)}_{\cal A}[[\zeta]]:=({\cal F}^{(3)}_{\cal A}[\zeta])[[\zeta]].$
\end{prop}
\begin{pf}
Note that by (\ref{mult}) we have the following identity
$$
\begin{array}{l}
{\rm exp}_{q^{2}}(-q^{-1}\zeta^{2}\mathop{}\hat{e}_{\alpha}
q^{-h_{\beta}^{\perp}}\otimes \hat{e}_{-\alpha-\beta})
=\exp_{q^{2}}(\zeta^{2}\mathop{}\hat{e}_{\alpha}
\displaystyle
\frac{H_{\alpha}^{\perp}-H_{-\beta}^{\perp}}{1+q}\otimes f_{0})\cdot\\[2ex]
\hphantom{aaaaaaaaaaaaaaaaaaaaaaaaa}
\exp_{q^{-2}}(-q^{-1}\zeta^{2}\mathop{}q^{h_{\alpha}^{\perp}}
\hat{e}_{\alpha}\otimes\hat{e}_{-\alpha-\beta})
\end{array}
$$
which allows to bring  $F_{3}^{\rm aff}$ to the following form
$$
\begin{array}{c}
(1\otimes 1-\zeta\mathop{}1\otimes f_{3}-\zeta^{2}\mathop{}
(h_{\beta}^{\perp}/2)_{q^{-2}}\otimes f_{2})_{q^{2}}
^{(-\frac 12 h_{\beta}^{\perp}\otimes 1)}\cdot\\[2ex]
(1\otimes 1+q^{-2}\zeta\mathop{}1\otimes\hat{e}_{-\alpha})_{q^{-2}}^{(-\frac 12 h_{\beta}^{\perp}\otimes 1)}\cdot
\exp_{q^{2}}(\zeta^{2}\mathop{}\hat{e}_{\alpha}\displaystyle
\frac{H_{\alpha}^{\perp}-H_{-\beta}^{\perp}}{1+q}\otimes f_{0})\cdot\\[2ex]
\exp_{q^{-2}}(-\zeta
\mathop{}q^{h_{\alpha}^{\perp}}\hat{e}_{\alpha}\otimes f_{1})\cdot
\exp_{q^{2}}(-q\zeta\mathop{}\hat{e}_{\alpha}\otimes f_{2})\cdot\\[2ex]
\hphantom{aaa}(1\otimes 1-q^{-2}\zeta\mathop{}
q^{h_{\alpha}^{\perp}}\otimes f_{3}
-\zeta^{2}\mathop(h_{\alpha}^{\perp}/2)_{q^{2}}\otimes f_{2})_{q^{-2}}
^{(\frac 12 h_{\alpha+\beta}^{\perp}\otimes 1)}
\end{array}
$$ 
correctly defined on ${\cal F}^{(3)}_{\cal A}[[\zeta]]$.
\end{pf}
The semi-classical twist $\overline{F^{\rm aff}_{3}}$ is obtained by specializing $q=1$: 
$$
\begin{array}{l}
(1\otimes 1-\zeta\mathop{}1\otimes\overline
{f^{\vphantom{f}}}_{3}-\frac 12\zeta^{2}\mathop{}
\overline{H^{\vphantom{f}}}_{\beta}^{\perp}\otimes\overline{f^{\vphantom{f}}}_{2})
^{(-\frac 12 \overline{H^{\vphantom{f}}}_{\beta}^{\perp}\otimes 1)}\cdot
(1\otimes 1+\zeta\mathop{} 1\otimes\overline{\vphantom{f}\hat{e}}_{-\alpha})^{(-\frac 12\overline{H^{\vphantom{X}}}_{\beta}^{\perp}\otimes 1)}\cdot\\[2ex]
\exp(\displaystyle\frac 12\mathop{}\zeta^{2}\mathop{}\overline{\vphantom{f}\hat{e}}_{\alpha}
(\overline{H^{\vphantom{f}}}_{\alpha}^{\perp}+\overline{H^{\vphantom{f}}}_{\beta}^{\perp})
\otimes\overline{f^{\vphantom{f}}}_{0})\cdot
\exp(-\zeta\mathop{}\overline{\vphantom{f}\hat{e}}_{\alpha}\otimes\overline{f^{\vphantom{f}}}_{1})\cdot
\exp(-\zeta\mathop{}\overline{\vphantom{f}\hat{e}}_{\alpha}\otimes\overline{f^{\vphantom{f}}}_{2})\cdot\\[2ex]
\hphantom{aaaaaaaaaaaaaaaaaaaaaa}
(1\otimes 1-\zeta\mathop{}1\otimes\overline{f^{\vphantom{f}}}_{3}-\frac 12\zeta^{2}\mathop{} 
\overline{H^{\vphantom{f}}}_{\alpha}^{\perp}\otimes\overline{f^{\vphantom{f}}}_{2})^{(\frac 12(\overline{H^{\vphantom{f}}}_{\beta}^{\perp}-\overline{H^{\vphantom{f}}}_{\alpha}^{\perp})\otimes 1)}.
\end{array}
$$
\section{Universal quantization of $n=3,4$ generalized Jordanian $r-$matrices}
Let us denote by $U^{\Psi_{n}}(\mathfrak{sl}_{n})[[\zeta]]$  
the completed universal enveloping algebra  $U(\mathfrak{sl}_{n})[[\zeta]]$ 
with the twisted coproduct  
$$
\Delta_{\Psi_{n}}(\cdot)=\Psi_{n}\Delta(\cdot)\Psi^{-1}_{n}.
$$
\begin{prop}
\label{prop3}
There exists  a twist $\Psi_{n}$ in $U(\mathfrak{sl}_{n})[[\zeta]]$ 
and a homomorphism $\iota_{n-1}$ such that 
$$
\iota_{n-1}:{\cal F}^{(n-1)}_{{\cal A},\mathop{}q=1}[[\zeta]]\rightarrow 
U^{\Psi_{n}}(\mathfrak{sl}_{n})[[\zeta]].
$$
\end{prop}
\begin{pf}
If $n=3$, then choose the extended Jordanian twist \cite{KLM}
$$
\Psi_{3}= {\rm exp}(\zeta\mathop{}E_{32}\otimes E_{13}e^{-\sigma_{12}^{-\zeta}})\cdot
{\rm exp}(D_{1}\otimes\sigma_{12}^{-\zeta})
$$
where $\sigma_{12}^{-\zeta}=\ln(1-\zeta E_{12})$ and $E_{ij}$ are the elements
of $U(\mathfrak{sl}_{n})$ corresponding to the 
elements of Cartan-Weyl basis of $\mathfrak{sl}_{n}$.
In the deformed $U^{\Psi_{3}}(\mathfrak{sl}_{3})[[\zeta]]$ we find the following elements 
and their coproducts
$$
\Delta_{\Psi_{3}}(E_{23})=E_{23}\otimes 1+1\otimes E_{23}-2\zeta\mathop{} 
D_{1}\otimes E_{13}e^{-\sigma_{12}^{-\zeta}}
$$
and $D_{1}$, $E_{13}e^{-\sigma_{12}^{-\zeta}}$ with primitive coproducts.
Define $\iota_{2}$ by its values on the generators of 
${\cal F}^{(2)}_{{\cal A},\mathop{}q=1}[[\zeta]]$:     
$$
\begin{array}{lclcl}
\iota_{2}(\overline{H^{\vphantom{f}}}_{\alpha})=-2 D_{1},&&
\iota_{2}(\overline{f^{\vphantom{f}}}_{0})=E_{13}e^{-\sigma_{12}^{-\eta}},&&
\iota_{2}(\overline{f^{\vphantom{f}}}_{1})=E_{23}
\end{array} 
$$ 
and extend it to ${\cal F}_{{\cal A},\mathop{}q=1}^{(2)}[[\zeta]]$ as a homomorphism 
into $U^{\Psi_{n}}(\mathfrak{sl}_{3})[[\zeta]]$. The defining relations (\ref{relations1})
and coproducts (\ref{corelations}) are preserved by $\iota_{2}$, thus the 
statement holds in this case.
Next, if $n=4$ then take 
$$
\Psi_{4}=\exp(\zeta\mathop{}E_{32}\otimes E_{13}e^{-\sigma_{12}^{-\zeta}}
+\zeta\mathop{}E_{42}\otimes E_{14}e^{-\sigma_{12}^{-\zeta}})\cdot
\exp(D_{1}\otimes\sigma_{12}^{-\zeta})
$$
and define $\iota_{3}:{\cal F}^{(3)}_{{\cal A},\mathop{}q=1}[[\zeta]]\rightarrow U^{\Psi_{4}}(\mathfrak{sl}_{4})[[\zeta]]$ by the following relations
$$
\begin{array}{rclcrclcrcl}
\iota_{3}(\overline{H^{\vphantom{f}}}_{\alpha})&=&D_{2}-2D_{3},&&
\iota_{3}(\overline{H^{\vphantom{f}}}_{\beta})&=&D_{3}-2 D_{2},&& 
\iota_{3}(\overline{f^{\vphantom{f}}}_{0})&=&-E_{14}e^{-\sigma_{12}^{-\zeta}}\\[2ex]
\iota_{3}(\overline{f^{\vphantom{f}}}_{1})&=&E_{24}^{\prime},&&
\iota_{3}(\overline{f^{\vphantom{f}}}_{2})&=&E_{13}e^{-\sigma_{12}^{-\zeta}}, &&
\iota_{3}(\overline{f^{\vphantom{f}}}_{3})&=&E_{23}\\[2ex]
\iota_{3}(\overline{\vphantom{f}\hat{e}}_{\alpha})&=&-E_{43}, &&
\iota_{3}(\overline{\vphantom{f}\hat{e}}_{-\alpha})&=&-E_{34}
\end{array}
$$
where $E^{\prime}_{24}=E_{24}-\zeta\mathop{}E_{34}E_{13}e^{-\sigma_{12}^{-\zeta}}$.
The only nonprimitive coproducts of the generators in $\iota_{3}({\cal F}^{(3)}_{{\cal A},\mathop{}q=1})[[\zeta]]$ 
are the following 
$$ 
\begin{array}{l}
\Delta_{\Psi_{4}}(E_{23})=\\[2ex]
\hphantom{aaa}E_{23}\otimes 1+1\otimes E_{23}+\zeta\mathop{}
(D_{3}-2D_{2})\otimes E_{13}e^{-\sigma_{12}^{-\zeta}}+
\zeta\mathop{}E_{43}\otimes E_{14}e^{-\sigma_{12}^{-\zeta}}\\[2ex]
\Delta_{\Psi_{4}}(E_{24}^{\prime})=\\[2ex]
\hphantom{aaa}E_{24}^{\prime}\otimes 1+1\otimes E_{24}^{\prime}
-\zeta\mathop{}(D_{2}+D_{3})\otimes E_{14}e^{-\sigma_{12}^{-\zeta}}-\zeta\mathop{}E_{13}e^{-\sigma_{12}^{-\zeta}}\otimes E_{34}.
\end{array}
$$
Let us consider the structure of ${\cal F}^{(3)}_{{\cal A},\mathop{}q=1}[[\zeta]]$. 
As a topological Hopf algebra it is the completion of 
${\cal F}^{(3)}_{{\cal A},\mathop{}q=1}[\zeta]:=$  
\begin{equation}
\{\sum_{l_{1},\ldots,l_{8}\ge 0}c_{l_{1},\ldots,l_{8}}(\overline{f^{\vphantom{f}}}_{2})^{l_{1}}(\overline{f^{\vphantom{f}}}_{0})^{l_{2}}
(\overline{\vphantom{f}\hat{e}}_{-\alpha})^{l_{3}}(
\overline{H^{\vphantom{f}}}_{\alpha}^{\perp})^{l_{4}}
(\overline{H^{\vphantom{f}}}_{\beta}^{\perp})^{l_{5}}
(\overline{\hat{e}\vphantom{f}}_{\alpha})^{l_{6}}
(\overline{f^{\vphantom{f}}}_{3})^{l_{7}}(\overline{f^{\vphantom{f}}}_{1})^{l_{8}}\}.
\label{basis123}
\end{equation}
The coproducts of the generators and the commutation relations are obtained from their quantum 
counterparts
$$
\begin{array}{c}
\Delta(\overline{f^{\vphantom{f}}}_{1})=
\overline{f^{\vphantom{f}}}_{1}\otimes 1+1\otimes\overline{f^{\vphantom{f}}}_{1}
-\zeta\mathop{}(\overline{H^{\vphantom{f}}}_{\alpha}+
\overline{H^{\vphantom{f}}}_{\beta})\otimes\overline{f^{\vphantom{f}}}_{0}+
\zeta\mathop{}
\overline{f^{\vphantom{f}}}_{2}\otimes
\overline{\vphantom{f}\hat{e}}_{-\alpha}\\[2ex]
\Delta(\overline{f^{\vphantom{f}}}_{3})=
\overline{f^{\vphantom{f}}}_{3}\otimes 1+1\otimes\overline{f^{\vphantom{f}}}_{3}
+\zeta\mathop{}\overline{H^{\vphantom{f}}}_{\beta}\otimes\overline{f^{\vphantom{f}}}_{2}
+\zeta\mathop{}\overline{\vphantom{f}\hat{e}}_{\alpha}\otimes 
\overline{f^{\vphantom{f}}}_{0},
\end{array}
$$ 
where we have written only the generators with nonprimitive coproducts,
$$
\begin{array}{lclcl}
[\overline{\vphantom{f}\hat{e}}_{\alpha},
\overline{\vphantom{f}\hat{e}}_{-\alpha}]=
\overline{H^{\vphantom{f}}}_{\alpha},&&[\overline{\vphantom{f}\hat{e}}_{\alpha},\overline{f^{\vphantom{f}}}_{0}]=
-\overline{f^{\vphantom{f}}}_{2}, && 
[\overline{\vphantom{f}\hat{e}}_{\alpha},\overline{f^{\vphantom{f}}}_{1}]=
\zeta\mathop{}\overline{f^{\vphantom{f}}}_{2}+\overline{f^{\vphantom{f}}}_{3}+\zeta\mathop{}
\overline{f^{\vphantom{f}}}_{2}\overline{H^{\vphantom{f}}}_{\alpha}
\\[2ex]
[\overline{\vphantom{f}\hat{e}}_{\alpha},\overline{f^{\vphantom{f}}}_{2}]=0,
&&[\overline{\vphantom{f}\hat{e}}_{\alpha},\overline{f^{\vphantom{f}}}_{3}]=0,&&
[\overline{f^{\vphantom{f}}}_{0},\overline{\vphantom{f}\hat{e}}_{-\alpha}]=
0\\[2ex]
[\overline{f^{\vphantom{f}}}_{1},\overline{\vphantom{f}\hat{e}}_{-\alpha}]=
\zeta\mathop{}
\overline{f^{\vphantom{f}}}_{0}\mathop{}\overline{\vphantom{f}\hat{e}}_{-\alpha},&&
[\overline{f^{\vphantom{f}}}_{2},\overline{\vphantom{f}\hat{e}}_{-\alpha}]=
\overline{f^{\vphantom{f}}}_{0},&&
[\overline{f^{\vphantom{f}}}_{3},\overline{\vphantom{f}\hat{e}}_{-\alpha}]
=-\zeta\overline{f^{\vphantom{f}}}_{0}-\overline{f^{\vphantom{f}}}_{1}
+\zeta\mathop{}\overline{f^{\vphantom{f}}}_{2}\mathop{}\overline{\vphantom{f}\hat{e}}_{-\alpha}\\[2ex]
[\overline{f^{\vphantom{f}}}_{0},\overline{f^{\vphantom{f}}}_{1}]=
-\zeta\mathop{}{\overline{f^{\vphantom{f}}}_{0}}^{2},&&
[\overline{f^{\vphantom{f}}}_{0}, \overline{f^{\vphantom{f}}}_{2}]=0,&&
[\overline{f^{\vphantom{f}}}_{0},\overline{f^{\vphantom{f}}}_{3}]=
\zeta\mathop{}\overline{f^{\vphantom{f}}}_{2}\overline{f^{\vphantom{f}}}_{0}\\[2ex]
[\overline{f^{\vphantom{f}}}_{1},\overline{f^{\vphantom{f}}}_{2}]=0,&&
[\overline{f^{\vphantom{f}}}_{1},\overline{f^{\vphantom{f}}}_{3}]=
\zeta\mathop{}\overline{f^{\vphantom{f}}}_{2}\overline{f^{\vphantom{f}}}_{1},&&
[\overline{f^{\vphantom{f}}}_{2},\overline{f^{\vphantom{f}}}_{3}]=
\zeta\mathop{}{\overline{f^{\vphantom{f}}}_{2}}^{2}.
\end{array}
$$
To see that any other relation in ${\cal F}^{(3)}_{{\cal A},\mathop{}q=1}[[\zeta]]$ 
follows from the introduced ones, we consider the quantum analogues 
of the commutation relations as the ordering rules in 
${\cal F}^{(3)}_{\cal A}[[\zeta]]$ with the  
normal ordering $f_{2}\prec f_{0}\prec
\hat{e}_{-\alpha}\prec H_{\alpha}
\prec H_{\beta}\prec
\hat{e}_{\alpha}\prec
f_{3}\prec f_{1}
$,
then by the Diamond lemma \cite{B} any monomial  
$m\in{\cal F}^{(3)}_{\cal A}[[\zeta]]$ can be 
brought to the form (\ref{basis321}), indeed otherwise there would exist a 
nontrivial relation between the ordered monomials and any such relation 
must have zero coefficients as the normal ordering 
in ${\cal F}^{(3)}_{\cal A}[[\zeta]]$ compatible with the
ordering in ${\cal D}^{(3)}_{\cal A}[[\zeta]]$:
$\hat{e}_{-\beta}\prec\hat{e}_{-\alpha-\beta}\prec\hat{e}_{-\alpha}\prec
H_{\alpha}\prec H_{\beta}\prec\hat{e}_{\alpha}\prec
\hat{e}_{\delta-\beta}\prec\hat{e}_{\delta-\alpha-\beta}$. Thus 
any monomial in ${\cal F}^{(3)}_{{\cal A},\mathop{}q=1}[[\zeta]]$
can be uniquely ordered as well and no additional relations 
come from ${\cal F}^{(3)}_{\cal A}[[\zeta]]$. Now it is direct to check 
that the extension of $\iota_{3}$ to a homomorphism preserves the 
coproducts and relations in ${\cal F}^{(3)}_{\cal A}[[\zeta]]$.
\end{pf}
The main new result of this paper is an explicit form of $F_{\mathfrak{p}}$ 
when $n=4$ is obtained as $F_{\mathfrak{p}}=(\iota_{3}\otimes\iota_{3})(\overline{F_{3}^{\rm aff}})\cdot\Psi_{4}$ 
and is the following
$$
\begin{array}{l}
(1\otimes 1-\zeta\mathop{}1\otimes E_{23}+\zeta^{2} D_{3}\otimes E_{13}e^{-\sigma_{12}^{-\zeta}})
^{(D_{3}\otimes 1)}\cdot
(1\otimes 1-\zeta\mathop{}1\otimes E_{34})^{(D_{3}\otimes 1)}\cdot\\[2ex]
\exp(-\zeta^{2}\mathop{}E_{43}(D_{2}+D_{3})\otimes E_{14}e^{-\sigma_{12}^{-\zeta}})\cdot
\exp(\zeta\mathop{}E_{43}\otimes(E_{24}-\zeta E_{34}E_{13}e^{-\sigma_{12}^{-\zeta}}))\cdot\\[2ex]
\exp(\zeta\mathop{}E_{43}\otimes E_{13}e^{-\sigma_{12}^{-\zeta}})\cdot
(1\otimes 1-\zeta\mathop{}1\otimes E_{23}+\zeta^{2}\mathop{}D_{2}\otimes E_{13}e^{-\sigma_{12}^{-\zeta}})^{((D_{2}-D_{3})\otimes 1)}\cdot\\[2ex]
\hphantom{aaaaaaaaaaaa}\exp(\zeta\mathop{}E_{32}\otimes E_{13}e^{-\sigma_{12}^{-\zeta}}+
\zeta\mathop{}E_{42}\otimes E_{14}e^{-\sigma_{12}^{-\zeta}})\cdot
\exp(D_{1}\otimes\sigma_{12}^{-\zeta}).
\end{array}
$$
Forming the universal $R-$matrix 
$R_{\mathfrak{p}}={F_{\mathfrak{p}}}_{21}F_{\mathfrak{p}}^{-1}$ we obtain 
in the first order in $\zeta$ the Gerstenhaber-Giaquinto $n=4$
$r_{\mathfrak{p}}-$matrix:
$$
D_{1}\wedge E_{12}+E_{14}\wedge E_{42}+E_{13}\wedge E_{32}
+D_{2}\wedge E_{23}+E_{24}\wedge E_{43}+D_{3}\wedge E_{34}
+E_{13}\wedge E_{43}.
\label{case4}  
$$
Note that there is a whole family of homomorphisms 
$\iota_{3}^{(a)}$:
$$
\begin{array}{rclcrclcrcl}
\iota_{3}^{(a)}(\overline{H^{\vphantom{f}}}_{\alpha})&=&D_{2}-2D_{3},&&
\iota_{3}^{(a)}(\overline{H^{\vphantom{f}}}_{\beta})&=&D_{3}-2 D_{2},&& 
\iota_{3}^{(a)}(\overline{f^{\vphantom{f}}}_{0})&=&\displaystyle\frac{1}{a}\mathop{}E_{14}e^{-\sigma_{12}^{-\zeta}}\\[2ex]
\iota_{3}^{(a)}(\overline{f^{\vphantom{f}}}_{1})&=&
-\displaystyle\frac{1}{a}\mathop{}E_{24}^{\prime},&&
\iota_{3}^{(a)}(\overline{f^{\vphantom{f}}}_{2})&=&E_{13}e^{-\sigma_{12}^{-\zeta}}, &&
\iota_{3}^{(a)}(\overline{f^{\vphantom{f}}}_{3})&=&E_{23}\\[2ex]
\iota_{3}^{(a)}(\overline{\vphantom{f}\hat{e}}_{\alpha})&=&a\mathop{}E_{43}, &&
\iota_{3}^{(a)}(\overline{\vphantom{f}\hat{e}}_{-\alpha})&=&\displaystyle\frac{1}{a}\mathop{}E_{34}
\end{array}
$$
and thus 
$(\iota_{3}^{(a)}\otimes\iota_{3}^{(a)})(F_{3}^{\rm aff})\cdot\Psi_{4}$  
for any  $a\ne 0$ leads to quantization of 
$$
D_{1}\wedge E_{12}+E_{13}\wedge E_{32}+E_{14}\wedge E_{42}
+D_{2}\wedge E_{23}+E_{24}\wedge E_{43}+\displaystyle\frac{1}{a}\mathop{}D_{3}\wedge E_{34}+
a\mathop{}E_{13}\wedge E_{43}.
$$ 
\begin{acknowledgements}
I would like to express my gratitude to V. Lyakhovsky, 
A. Stolin and V. Tolstoy for the valuable discussions at the different stages of the work.
\end{acknowledgements}

\end{article}
\end{document}